\documentclass[12pt, reqno]{amsart}
\usepackage{amssymb}
\usepackage{amscd}
\usepackage{verbatim,ifthen}
\usepackage{color}
\usepackage{latexsym}
\usepackage{tikz}
\usepackage{tikz-cd}
\usepackage{mathrsfs}
\usepackage{wrapfig}
\usetikzlibrary{shapes}
\usepackage{color}
\usetikzlibrary{arrows.meta}
\usepackage{bbm}
\usetikzlibrary{matrix}
\usetikzlibrary{calc}
\usetikzlibrary{arrows,intersections}
\usepackage{pgfplots}
\usepackage{multicol}
\usepackage{array}
\newcolumntype{M}[1]{>{\centering\arraybackslash}m{#1}}

\usepackage[colorlinks, linkcolor=black, citecolor=blue, linktocpage,urlcolor=blue,backref=page]{hyperref}
\renewcommand*{\backref}[1]{}
\renewcommand*{\backrefalt}[4]{[{\tiny%
    \ifcase #1 Not cited.%
          \or Cited on page~#2.%
          \else Cited on pages #2.%
    \fi%
    }]}
\addtocontents{toc}{\setcounter{tocdepth}{1}}

\usepackage{amsmath}

\numberwithin{equation}{section}

\let\oldtocsection=\tocsection
 
\let\oldtocsubsection=\tocsubsection

\renewcommand{\tocsection}[2]{\hspace{0em}\oldtocsection{#1}{#2}}
\renewcommand{\tocsubsection}[2]{\hspace{1em}\oldtocsubsection{#1}{#2}}

\def\XXint#1#2#3{{\setbox0=\hbox{$#1{#2#3}{\int}$ }
\vcenter{\hbox{$#2#3$ }}\kern-.6\wd0}}

\usepackage{eucal}
\usepackage{calc}  
\usepackage{enumitem} 
\usepackage{tensor}
\usepackage{graphicx,wrapfig,lipsum}
\usepackage{etoolbox}
\usepackage{marginnote}
\usepackage{lipsum}
\makeatletter
\patchcmd{\@mn@margintest}{\@tempswafalse}{\@tempswatrue}{}{}
\patchcmd{\@mn@margintest}{\@tempswafalse}{\@tempswatrue}{}{}
\reversemarginpar 
\makeatother
\usepackage{scrextend}

\makeatletter
\DeclareRobustCommand\widecheck[1]{{\mathpalette\@widecheck{#1}}}
\def\@widecheck#1#2{%
    \setbox\z@\hbox{\m@th$#1#2$}%
    \setbox\tw@\hbox{\m@th$#1%
       \widehat{%
          \vrule\@width\z@\@height\ht\z@
          \vrule\@height\z@\@width\wd\z@}$}%
    \dp\tw@-\ht\z@
    \@tempdima\ht\z@ \advance\@tempdima2\ht\tw@ \divide\@tempdima\thr@@
    \setbox\tw@\hbox{%
       \raise\@tempdima\hbox{\scalebox{1}[-1]{\lower\@tempdima\box
\tw@}}}%
    {\ooalign{\box\tw@ \cr \box\z@}}}
\makeatother
\title{$(\varepsilon,\delta)$--Quasi-Negative Curvature and Positivity of the Canonical Bundle}
\author{Kyle Broder}
\address{The University of Queensland,  St. Lucia,  QLD 4067, Australia}
\email{k.broder@uq.edu.au}
\author{Kai Tang}
\address{College of Mathematics and Computer Science, Zhejiang Normal University, Jinhua, Zhejiang, 321004, China}
\email{kaitang001@zjnu.edu.cn}
\thanks{The first named author was supported by the Australian Government through the Australian Research Council's Discovery Projects funding scheme (project DP220102530). 
The second named author was supported by National Natural Science
Foundation of China (Grant No.12001490).
}
\begin{document}

\maketitle

\begin{abstract}
A recent theorem of Diverio--Trapani and Wu--Yau asserts that a compact K\"ahler manifold with a K\"ahler metric of quasi-negative holomorphic sectional curvature is projective and canonically polarized.  This confirms a long-standing conjecture of Yau. We consider the notion of $(\varepsilon,\delta)$--quasi-negativity, generalizing quasi-negativity, and obtain gap-type theorems for $\int_X c_1(K_X)^n>0$ in terms of the real bisectional curvature and weighted orthogonal Ricci curvature. These theorems are also a generalization of
that results in \cite{ZhangZheng} by Zhang-Zheng and in \cite{ChuLeeTam} by Chu-Lee-Tam.
\end{abstract}

\section{Introduction}
The Wu--Yau theorem \cite{Campana, Wong, HeierLuWong, WuYau1, WuYau2, TosattiYang} asserts that a compact K\"ahler manifold with negative holomorphic sectional curvature is projective and canonically polarized.  In particular, the canonical bundle $K_X$ is ample and there is a unique invariant K\"ahler--Einstein metric \cite{Aubin,Yau1976}.  Diverio--Trapani \cite{DiverioTrapani} relaxed the strict negativity of the holomorphic sectional curvature to quasi-negativity (non-positive everywhere and negative at one point). The crux of the argument is the positivity of the top intersection number \begin{eqnarray}\label{PositiveIntersectionNumber}
\int_X c_1(K_X)^n \ > \ 0.
\end{eqnarray} Indeed,  from \cite[Theorem 0.5]{DemaillyPaun}, \eqref{PositiveIntersectionNumber} implies that the canonical bundle $K_X$ is big if we assume that $K_X$ is nef. Hence, $X$ is Moishezon, and by Moishezon's theorem \cite{Moishezon},  $X$ is projective and of general type.  A projective manifold of general type with no rational curves has ample canonical bundle (see \cite[Lemma 2.1]{DiverioTrapani} or \cite[p. 219]{Debarre}).

In their recent work, Zhang-Zheng \cite{ZhangZheng} considered a natural notion of almost quasi-negative holomorphic sectional curvature and extended these theorems to compact K\"ahler manifolds of almost
quasi-negative holomorphic sectional curvature. They also obtained a gap-type theorem for
the inequality $\int_X c_1(K_X)^n  > 0$ in terms of the holomorphic sectional curvature on compact K\"ahler manifolds.
We will extend these results of Diverio--Trapani \cite{DiverioTrapani} and  Zhang-Zheng \cite{ZhangZheng} to the Hermitian category.  We remind the reader that for non-K\"ahler Hermitian metrics,  the holomorphic sectional curvature is naturally replaced by the \textit{real bisectional curvature} $$\text{RBC}_{\omega}(\xi) \ : = \ \frac{1}{| \xi |^2} \sum_{i,j,k,\ell} R_{i \overline{j} k \overline{\ell}} \xi^{i \overline{j}} \xi^{k \overline{\ell}},$$ where $\xi$ is a non-negative $(1,1)$--tensor.  The real bisectional curvature was introduced by Yang--Zheng \cite{YangZhengRBC} in a frame-dependent manner. The above definition is also due to Lee--Streets \cite{LeeStreets}, where it is referred to as the `complex curvature operator'.

When the metric is K\"ahler (or, more generally,  K\"ahler-like \cite{YangZhengCurvature}), the real bisectional curvature is equivalent to the altered holomorphic sectional curvature formally defined in \cite{BroderTangAltered}, which is comparable to the familiar holomorphic sectional curvature in the sense that they always have the same sign. In general, however, the real bisectional curvature is stronger than the holomorphic sectional curvature.  It is not strong enough, however, to dominate the Ricci curvatures (see \cite{YangZhengRBC}). For a discussion of these curvatures for connections more general than the Chern connection,  namely,  the Gauduchon connections, we invite the reader to see \cite{BroderStanfield}.

To state the main results, we introduce the following terminology (building from \cite{ZhangZheng}): For a positive constant $\delta>0$ and a non-empty open set $\mathcal{U} \subset X$ we say that a Hermitian metric $\eta$ has \textit{$(\varepsilon,\delta)$--quasi-negative real bisectional curvature} (relative to $\mathcal{U}$) if there is a sufficiently small $\varepsilon>0$ such that $\text{RBC}_{\eta} \leq \varepsilon$ on $X$ and if $\text{RBC}_{\eta} \leq - \delta$ on $\mathcal{U}$.   For positive constants $\delta_1, \delta_2 >0$, we say that a Hermitian metric $\eta$ on a compact K\"ahler manifold $(X, \omega)$ is (i) \textit{$\delta_1$--bounded} if there is a smooth function $\psi : X \to \mathbf{R}$ such that $\eta \ \leq \ \delta_1 \omega + dd^c \psi$, and (ii) \textit{$\delta_2$--volume non-collapsed} on an open set $\mathcal{U} \subset X$ if $\eta^n \geq \delta_2 \omega^n$. If a Hermitian metric is both $\delta_1$--bounded and $\delta_2$--volume non-collapsed, we say that $\eta$ has $(\delta_1, \delta_2)$--bounded geometry.

The first main theorem of the present manuscript is the following:

\subsection*{Theorem 1.1}\label{Theorem2}
Let $X$ be a compact K\"ahler manifold.  Suppose there is a Hermitian metric $\eta$ with $(\delta_1, \delta_2)$--bounded geometry and $(\varepsilon, \delta)$--quasi-negative real bisectional curvature.  Then $$\int_X c_1(K_X)^n \ > \ 0$$ If one further assumes that  $K_X$ is nef, then  $K_X$  is big and $X$ is projective. \\

The above theorem is a result in this direction that does not require the curvature of the Hermitian metric to have a sign.  There are also many interesting results in \cite{Zhang,Zhang1,ZhangZheng,TangkRicci} without assuming a pointwise signed curvature condition. Further,  let us emphasize that although the $\delta_1$--boundedness property is pointwise, it can be interpreted as a bound on the Bott--Chern class represented by the metric, not at the level of the metric.  On the other hand, $\delta_2$--volume non-collapsing is a local weighted volume non-collapsing condition; for instance, it holds if $\eta^n \geq \delta_2 \omega_0^n$ on $\mathcal{U}$.\\

Given the perplexing relationship between the holomorphic sectional curvature and the Ricci curvature (c.f., \cite{BroderRemarks}),  one can attempt to interpolate between these curvatures by considering the following curvature constraint: In \cite{BT} (see also \cite{ChuLeeTam}), the  authors introduced, for $\alpha, \beta \in \mathbf{R}$ the \textit{weighted orthogonal Ricci curvature} $$\text{Ric}_{\alpha, \beta}^{\perp}(v) \ : = \ \frac{\alpha}{| v |_{\omega_0}^2} \text{Ric}_{\omega_0}(v, \overline{v}) + \beta \text{HSC}_{\omega_0}(v),$$ for a $(1,0)$--tangent vector $v \in T^{1,0}X$.

From Siu's resolution \cite{SiuMoishezon} of the Grauert--Riemenschneider conjecture, a compact K\"ahler manifold with quasi-negative Ricci curvature has big canonical bundle.  On the other hand, the Diverio--Trapani \cite{DiverioTrapani} extension of the Wu--Yau theorem shows that a compact K\"ahler manifold with quasi-negative holomorphic sectional curvature has big canonical bundle.  The following theorem shows that these results can be interpolated via the weighted orthogonal Ricci curvature, and more generally, extended to the $(\varepsilon,\delta)$--quasi-negative situation:

\subsection*{Theorem 1.2}\label{Theorem3}
Let $(X, \omega)$ be a compact K\"ahler manifold with $(\varepsilon, \delta)$--quasi-negative weighted orthogonal Ricci curvature for $\alpha, \beta \in \mathbf{R}_{\geq 0}$ with at least one positive. Then \begin{eqnarray*}
\int_X c_1(K_X)^n & > & 0
\end{eqnarray*}
 If one further assumes that  $K_X$ is nef, then  $K_X$  is big and $X$ is projective. \\

\nameref{Theorem3} also extends a result of Chu--Lee--Tam \cite{ChuLeeTam}, where they showed that if a compact K\"ahler manifold $X$ has a K\"ahler metric of quasi-negative weighted orthogonal Ricci curvature (for $\alpha, \beta \in \mathbf{R}_{>0}$), then $\int_X c_1(K_X)^n>0$.  \\

An alternative means of interpolating between the holomorphic sectional curvature and the Ricci curvature was described by Ni \cite{NikRicci}: The so-called \textit{$k$--Ricci curvature} of a K\"ahler metric is the Ricci curvature of the $k$--dimensional holomorphic subspaces of $T^{1,0}X$.  By making use of a result in \cite{ChuLeeTam}, the method used to prove \nameref{Theorem2} and \nameref{Theorem3} also provides the following (see also \cite{TangkRicci}):

\subsection*{Theorem 1.3}\label{Theorem4}
Let $(X, \omega)$ be a compact K\"ahler manifold with $(\varepsilon,\delta)$--quasi-negative $k$--Ricci curvature. Then $$\int_X c_1(K_X)^n \ > \ 0$$ If one further assumes that  $K_X$ is nef, then  $K_X$  is big and $X$ is projective.

\subsection*{Acknowledgements}
The first named author would like to thank his former Ph.D.  advisors Ben Andrews and Gang Tian for their support and encouragement.  He would also like to thank Simone Diverio,  Jeffrey Streets,  Finnur L\'arusson, Ramiro Lafuente,  and James Stanfield for many valuable discussions and communications. The authors would like to thank Yashan Zhang for his support, and for many invaluable discussions and suggestions.

\section{Proof of Theorem 1.1}
To show that $\int_X c_1(K_X)^n >0$,  it suffices to obtain the estimate \begin{eqnarray}\label{WTS}
\int_X (-\text{Ric}_{\omega_0})^n = (2\pi)^n \int_X c_1(K_X)^n >0.
\end{eqnarray} 

Let $\rho : \mathbf{R} \to \mathbf{R}$ be the function defined by $\rho(t) = \frac{1}{n}$ for $t \leq 0$ and $\rho(t) = 1$ for $t > 0$.  For the Hermitian metric $\eta$ in the statement of \nameref{Theorem2},  let $\kappa_{\eta} : X \to \mathbf{R}$ be the function \begin{eqnarray*}
\kappa_{\eta}(x) \ : = \ \rho \left( \max_{(\vartheta,v_x) \in \mathcal{F}_X \times \mathbf{R}^n} \text{RBC}_{\eta}(\vartheta,v_x) \right) \cdot \max_{(\vartheta,v_x) \in \mathcal{F}_X \times \mathbf{R}^n} \text{RBC}_{\eta}(\vartheta,v_x),
\end{eqnarray*}

where $\vartheta$ is a unitary frame (i.e., a section of the unitary frame bundle $\mathcal{F}_X$ and $v_x\in \mathbf{R}^n$).

Let us also introduce the notation \begin{eqnarray}\label{MU}
\mu_{\eta} & : = & \max_{x \in X} \max_{(\vartheta,v_x) \in \mathcal{F}_X \times \mathbf{R}^n} \text{RBC}_{\eta}(\vartheta,v_x).
\end{eqnarray}

To establish \eqref{WTS},  we will show that there are constants $\varepsilon, c_3, c_4 >0$ such that \begin{eqnarray}
\int_X (-\text{Ric}_{\omega_0})^n & \geq & \int_X (n \delta_1 \kappa_{\eta} \omega_0 - \text{Ric}_{\omega_0})^n - c_4 \varepsilon \label{WTS1} \\
& \geq & c_3 - c_4 \varepsilon, \label{WTS2}
\end{eqnarray}

where $\varepsilon >0$ can be chosen such that $c_3 - c_4 \varepsilon >0$.  \\

To this end, consider the twisted Wu--Yau continuity method, given by the complex Monge--Amp\`ere equation \begin{eqnarray}\label{WYCM}
(t (\omega_0 + \delta_1^{-1} dd^c  \psi) - \text{Ric}_{\omega_0} + dd^c  \varphi_t)^n &=& e^{\varphi_t} \omega_0^n.
\end{eqnarray}

From the assumption of $\delta_1$--boundedness, we see that \begin{eqnarray}\label{Delta1}
\omega_0 + \delta_1^{-1} dd^c \psi & \geq & \delta_1^{-1} \eta.
\end{eqnarray} Set $\omega_t : = t (\omega_0 + \delta_1^{-1} dd^c  \psi) - \text{Ric}_{\omega_0} + dd^c \varphi_t$,  allowing us to write \eqref{WYCM} as \begin{eqnarray*}\label{WYCM2}
\omega_t^n &=& e^{\varphi_t} \omega_0^n.
\end{eqnarray*}

From \eqref{Delta1},  the metrics $\omega_t$ afford the lower bound \begin{eqnarray*}
\text{Ric}_{\omega_t} & \geq & - \omega_t + t \delta_1^{-1} \eta.
\end{eqnarray*}

We will make use of the following Schwarz lemma for holomorphic maps between Hermitian manifolds due to Yang--Zheng \cite{YangZhengRBC} (see also \cite{Royden, BroderSBC,BroderSBC2}):

\subsection*{Lemma 2.1}\label{SchwarzLemma}
(Yang--Zheng). Let $f : (X, \omega_g) \to (Y, \eta)$ be a holomorphic map from a K\"ahler manifold to a Hermitian manifold.  If there are constants $C_1, C_2 \in \mathbf{R}$ such that $\text{Ric}_{\omega_g} \geq - C_1  \omega + C_2 f^{\ast} \omega_h$ and $\text{RBC}_{\eta} \leq \mu_{\eta}$, then  \begin{eqnarray*}
\Delta_{\omega_g} \log \text{tr}_{\omega_g}(f^{\ast} \eta) & \geq & - C_1 +  \left( - \mu_{\eta} + \frac{C_2}{n} \right) \text{tr}_{\omega_g}(f^{\ast} \eta).
\end{eqnarray*} 

\hfill

Taking $f$ to be the identity map,  $C_1 =1$,  and $C_2 = t \delta_1^{-1}$ in \nameref{SchwarzLemma}, we have \begin{eqnarray}\label{CLEstimate}
\Delta_{\omega_t} \log \text{tr}_{\omega_t}(\eta) & \geq & (-\mu_{\eta} + t(n \delta_1)^{-1} ) \text{tr}_{\omega_t}(\eta) - 1.
\end{eqnarray}

By the maximum principle, \begin{eqnarray}\label{C2Estimate}
\sup_X \text{tr}_{\omega_t}(\eta) & \leq & \frac{n\delta_1}{t-\mu_{\eta} n \delta_1}.
\end{eqnarray}

For $t > n \delta_1 \mu_{\eta}$, the estimate \eqref{C2Estimate} is independent of $t$.  As a consequence, the continuity method admits a smooth solution for $t > n \delta_1 \mu_{\eta}$ (c.f., \cite{TosattiYang, Zhang}).  Introduce the potential $u_t : = \varphi_t + t \delta_1^{-1} \psi$.  The crux of the argument is to estimate $\sup_{\mathcal{U}} u_t$ from below, and $\sup_X u_t$ from above.  Indeed,  the constant $c_3$ in \eqref{WTS2} is given by \begin{eqnarray*}
c_3 & : = & \liminf_{t \to n \delta_1 \mu_{\eta}} \int_X e^{u_t} \omega_0^n,
\end{eqnarray*}

where $\mu_{\eta}$ is defined in \eqref{MU}. The positivity of $c_3$ demands $u_t$ to not be identically $-\infty$, while the finiteness of $u_t$ requires an upper bound on $\sup_X u_t$.  Before obtaining the upper bound on $\sup_X u_t$, we first note that in contrast to the situations considered by Wu--Yau \cite{WuYau1}, Tosatti--Yang \cite{TosattiYang}, and Diverio--Trapani \cite{DiverioTrapani},  we do not have a metric with signed curvature.  To handle the positive curvature contributions,  we make use of \textit{Tian's $\alpha$--invariant} \cite{TianAlpha} (see also \cite{ZhangZheng}):

\subsection*{Proposition 2.2}
(Tian's $\alpha$--invariant). Let $(X, \omega)$ be a compact K\"ahler manifold. There is a positive number $\alpha = \alpha(X, [\omega])$ depending only on $X$ and the K\"ahler class $[\omega]$ such that for any $\beta \in (0,\alpha)$ and any $u \in \mathcal{C}^{\infty}(X, \mathbf{R})$ with $\omega + dd^c u > 0$ and $\sup_X u =0$,  there is a constant $C_{\beta} = C_{\beta}(\omega)$ such that $$\int_X e^{-\beta u } \omega^n \ \leq \ C_{\beta}.$$

\hfill

We now obtain an upper bound on $\sup_X u_t$:

\subsection*{Lemma 2.3}
For $t \in (n \delta_1 \mu_{\eta}, 2n \delta_1 \mu_{\eta}]$, we have \begin{eqnarray}\label{suput}
\sup_{x \in X} u_t(x) & \leq & \log(2n \delta_1 \varepsilon + b_0)^n \ \leq \  \log(c_0 + b_0)^n,
\end{eqnarray}

where $c_0>0$ is such that Tian's $\alpha$--invariant satisfies $\alpha(X, c_0\omega_0) \geq 2$, and $\varepsilon \leq c_0/(2n\delta_1)$.

\begin{proof}
Since $u_t = \varphi_t  + t \delta_1^{-1} \psi$,  \eqref{WYCM2} reads\begin{eqnarray}\label{WYCM3}
e^{u_t - t \delta_1^{-1} \psi} \omega_0^n \ = \ e^{\varphi_t} \omega_0^n \ = \ (t \omega_0 - \text{Ric}_{\omega_0} + dd^c  u_t)^n.
\end{eqnarray}

Since the metric $\omega_0$ is smooth, and $X$ is compact, there is some $b_0 \geq 0$ such that \begin{eqnarray}\label{b0}
\text{Ric}_{\omega_0} & \geq & - b_0 \omega_0.
\end{eqnarray}

From \eqref{WYCM3}, \eqref{b0},  and the fact that $\psi \leq 0$,  at the point $x \in X$ where $u_t$ achieves its maximum,  we have \begin{eqnarray*}
\sup_{x \in X} u_t(x) & \leq & \sup_X \left( t \delta_1^{-1} \psi + \log \frac{(t \omega_0 - \text{Ric}_{\omega_0} + dd^c u_t)^n}{\omega_0^n} \right) \  \leq  \ \log(t_0 + b_0)^n.
\end{eqnarray*}

Since $t  \leq 2n \delta_1 \mu_{\eta} \leq 2n \delta_1 \varepsilon$,  this proves \eqref{suput}.

\end{proof}

We now want to estimate $\sup_{\mathcal{U}} u_t$ from below.  

\subsection*{Lemma 2.4}
\begin{eqnarray}\label{LOGLOG}
\sup_{\mathcal{U}} u_t & \geq & n \log(n) + \log \left(  \frac{- \int_{\mathcal{U}}  \kappa_{\eta} e^{t(n\delta_1)^{-1} \psi} \left( \frac{\eta^n}{\omega_0^n} \right)^{\frac{1}{n}} \omega_t^n }{\int_X \omega_t^n} \right).
\end{eqnarray}

\begin{proof}
Start by integrating \eqref{CLEstimate} over $X$ with respect to the volume form $\omega_t^n$. By the divergence theorem,   \begin{eqnarray*}
\int_X \omega_t^n & \geq & \int_X (-\kappa_{\eta} + t(n\delta_1)^{-1}) \text{tr}_{\omega_t}(\eta) \omega_t^n \\
& \geq & \int_{\mathcal{U}} (- \kappa_{\eta} + t (n \delta_1)^{-1}) \text{tr}_{\omega_t}(\eta) \omega_t^n \ \geq \  - \int_{\mathcal{U}} \kappa_{\eta} \text{tr}_{\omega_t}(\eta) \omega_t^n,
\end{eqnarray*}

using the fact that $(-\kappa_{\eta} + t(n\delta_1)^{-1})\text{tr}_{\omega_t}(\eta) >0$.  By the arithmetic-geometric mean inequality,  \begin{eqnarray}
-\int_{\mathcal{U}} \kappa_{\eta} \text{tr}_{\omega_t}(\eta) \omega_t^n & \geq & - n \int_{\mathcal{U}} \kappa_{\eta} \left( \frac{\eta^n}{\omega_t^n} \right)^{\frac{1}{n}} \omega_t^n \nonumber \\
&=& -n \int_{\mathcal{U}} \kappa_{\eta} \left( \frac{\omega_0^n}{\omega_t^n} \right)^{\frac{1}{n}} \left( \frac{\eta^n}{\omega_0^n} \right)^{\frac{1}{n}} \omega_t^n. \label{AMGM}
\end{eqnarray}

Using \eqref{WYCM3} in \eqref{AMGM}, we achieve the estimate \begin{eqnarray}
\int_X \omega_t^n & \geq & -n \int_{\mathcal{U}} \kappa_{\eta} e^{-\frac{1}{n}(u_t - t \delta_1^{-1} \psi)} \left( \frac{\eta^n}{\omega_0^n} \right)^{\frac{1}{n}} \omega_t^n \nonumber \\
& \geq & - n e^{-\frac{1}{n} \sup_{\mathcal{U}} u_t } \int_{\mathcal{U}} \kappa_{\eta} e^{t(n \delta_1)^{-1} \psi} \left( \frac{\eta^n}{\omega_0^n} \right)^{\frac{1}{n}} \omega_t^n. \label{LBVol}
\end{eqnarray}

Then \eqref{LOGLOG} follows from \eqref{LBVol}.
\end{proof}

From the above lemma, it suffices to estimate \begin{eqnarray}\label{RATIO}
\frac{- \int_{\mathcal{U}}  \kappa_{\eta} e^{t(n\delta_1)^{-1} \psi} \left( \frac{\eta^n}{\omega_0^n} \right)^{\frac{1}{n}} \omega_t^n }{\int_X \omega_t^n} 
\end{eqnarray}

from below.  The following lemma gives an estimate for the numerator: 

\subsection*{Lemma 2.5}
There are positive constants $c_1, \delta_2>0$ such that 
\begin{eqnarray}\label{NUMERATOR}
- \int_{\mathcal{U}}  \kappa_{\eta} e^{t(n\delta_1)^{-1} \psi} \left( \frac{\eta^n}{\omega_0^n} \right)^{\frac{1}{n}} \omega_t^n  & \geq & \delta_2 \delta \int_{\mathcal{U}} e^{u_t^{\ast}} \omega_0^n \ \geq \ \frac{c_1 \delta_2^{\frac{1}{n}}\delta}{n}.
\end{eqnarray}

\begin{proof}
Let us write $u_t^{\ast} : = u_t - \sup_X u_t$, so that $\sup_X u_t^{\ast} \leq 0$. In this notation, \eqref{RATIO} reads \begin{eqnarray}
\frac{- \int_{\mathcal{U}}  \kappa_{\eta} e^{t(n\delta_1)^{-1} \psi} \left( \frac{\eta^n}{\omega_0^n} \right)^{\frac{1}{n}} \omega_t^n }{\int_X \omega_t^n}  &=& \frac{- \int_{\mathcal{U}} \kappa_{\eta} e^{u_t} e^{-\left( 1 - \frac{1}{n} \right) t \delta_1^{-1} \psi} \left( \frac{\eta^n}{\omega_0^n} \right)^{\frac{1}{n}} \omega_0^n}{\int_X e^{u_t} e^{-t \delta_1^{-1} \psi} \omega_0^n} \nonumber \\
&=& \frac{-\int_{\mathcal{U}} \kappa_{\eta} e^{u_t^{\ast}} e^{- \left( 1 - \frac{1}{n} \right) t \delta_1^{-1} \psi}\left( \frac{\eta^n}{\omega_0^n} \right)^{\frac{1}{n}} \omega_0^n}{\int_X e^{u_t^{\ast}} e^{-t \delta_1^{-1} \psi} \omega_0^n} \label{FRAC}
\end{eqnarray}

Since $\psi \leq 0$,  and $t \in (n \delta_1 \mu_{\eta}, 2n \delta_1 \mu_{\eta}]$, we have \begin{eqnarray}
- \int_{\mathcal{U}} \kappa_{\eta} e^{u_t^{\ast}} e^{-\left( 1 - \frac{1}{n} \right) t \delta_1^{-1} \psi} \left( \frac{\eta^n}{\omega_0^n} \right)^{\frac{1}{n}} \omega_0^n & \geq &  - \int_{\mathcal{U}} \kappa_{\eta} e^{u_t^{\ast}} e^{(1-n) \kappa_{\eta} \psi} \left( \frac{\eta^n}{\omega_0^n} \right)^{\frac{1}{n}} \omega_0^n \nonumber \\
& \geq &  \delta \int_{\mathcal{U}} e^{u_t^{\ast}} e^{(1-n)\kappa_{\eta} \psi} \left( \frac{\eta^n}{\omega_0^n} \right)^{\frac{1}{n}} \omega_0^n, \label{DELTA2}
\end{eqnarray}

where the last inequality follows from the negative curvature estimate $\text{RBC}_{\eta} \leq - \delta_3$ on $\mathcal{U}$.  Let \begin{eqnarray}\label{C1}
c_1 &: =& \inf \left \{ \int_{\mathcal{U}} e^v \omega_0^n \ : \ v \in \text{PSH}_{(c_0+b_0)\omega_0}(X), \ \sup_X v =0 \right \}.
\end{eqnarray} Then $c_1$ is a positive constant depending only on $\mathcal{U}$ and $\omega_0$.  Since $\eta$ has $(\delta_1, \delta_2)$--bounded geometry, we have \begin{eqnarray*}
-\int_{\mathcal{U}} \kappa_{\eta} e^{u_t^{\ast}} e^{-(n-1) \mu_{\eta} \psi} \left( \frac{\eta^n}{\omega_0^n} \right)^{\frac{1}{n}} \omega_0^n & \geq & \frac{\delta_2^{\frac{1}{n}} \delta}{n} \int_{\mathcal{U}} e^{u_t^{\ast}} \omega_0^n \ \geq \ \frac{c_1 \delta_2^{\frac{1}{n}}\delta}{n}.
\end{eqnarray*}

From \eqref{DELTA2} and \eqref{C1}, this gives the desired lower bound for the numerator in \eqref{FRAC}.
\end{proof}

We now complete the estimate for lower bound on $\sup_{\mathcal{U}} u_t$: 

\subsection*{Lemma 2.6}
There are constants $c_1, c_2, \delta_1, \delta_2 >0$ such that 
\begin{eqnarray*}
\sup_{\mathcal{U}} u_t & \geq & n \log(n) + \log \left( \frac{c_1 \delta_2 \delta}{c_2} \right).
\end{eqnarray*}

\begin{proof}
For the denominator in \eqref{FRAC},  we will again use Tian's $\alpha$--invariant \cite{TianAlpha}. Indeed,  first observe that, since  $0 < t \leq 2n \delta_1 \mu_{\eta} \leq c_0$,  \eqref{Delta1} implies \begin{eqnarray*}
c_0 \omega_0 + t \delta_1^{-1} dd^c \psi & \geq & t \omega_0 + t \delta_1^{-1} dd^c \psi \ \geq \ t \delta_1^{-1} \eta  \ > \ 0.
\end{eqnarray*}

In other words, $t \delta_1^{-1} \psi$ is $c_0 \omega_0$--plurisubharmonic.  Hence, since $c_0>0$ is chosen such that Tian's $\alpha$--invariant satisfies $\alpha(X, c_0 \omega_0) \geq 2$,  we have \begin{eqnarray}\label{DENOMINATOR}
\int_X e^{u_t^{\ast}} e^{-t \delta_1^{-1} \psi} \omega_0^n & \leq & \int_X e^{-t \delta_1^{-1} \psi} \omega_0^n  \ = \ c_0^{-n} \int_X e^{-t \delta_1^{-1} \psi} (c_0 \omega_0)^n  \ \leq \ c_2,
\end{eqnarray}

for some constant $c_2>0$ depending only on $c_0$ and $\omega_0$. Combining \eqref{RATIO}, \eqref{FRAC}, \eqref{NUMERATOR}, and \eqref{DENOMINATOR}, we see that \begin{eqnarray}\label{FRACESTIMATE}
\frac{- \int_{\mathcal{U}}  \kappa_{\eta} e^{t(n\delta_1)^{-1} \psi} \left( \frac{\eta^n}{\omega_0^n} \right)^{\frac{1}{n}} \omega_t^n }{\int_X \omega_t^n}  & \geq & \frac{c_1 \delta_2^{\frac{1}{n}} \delta}{n c_2}.
\end{eqnarray}

Inserting \eqref{FRACESTIMATE} into \eqref{LOGLOG}, \begin{eqnarray*}
\sup_{\mathcal{U}} u_t & \geq & n \log(n) + \log \left( \frac{c_1 \delta_2^{\frac{1}{n}} \delta}{n c_2} \right).
\end{eqnarray*}
\end{proof}

We now complete the proof: Let $c_4>0$ be the constant such that \begin{eqnarray*}
\int_X \sum_{k=1}^n \binom{n}{k} (2n \delta_1 \varepsilon \omega_0)^k \wedge (-\text{Ric}_{\omega_0})^{n-k} & \leq &c_4 \varepsilon.
\end{eqnarray*}

Then \begin{eqnarray*}
\int_X (2\pi c_1(K_X))^n \ =\  \int_X (-\text{Ric}_{\omega_0})^n & \geq & \int_X (2n \delta_1 \varepsilon \omega_0 - \text{Ric}_{\omega_0})^n -c_4 \varepsilon \\
& \geq & \lim_{t \searrow n \delta_1 \mu_{\eta}} (t \omega_0 - \text{Ric}_{\omega_0} + dd^c u_t)^n - c_4 \varepsilon \\
&=& \lim_{t \searrow n \delta_1 \mu_{\eta}} \int_X e^{u_t} e^{-t \delta_1^{-1} \psi} \omega_0^n - c_4 \varepsilon \\
& \geq & \liminf_{t \searrow n \delta_1 \mu_{\eta}} \int_X e^{u_t} \omega_0^n - c_4 \varepsilon  \\
&=& c_3 - c_4 \varepsilon.
\end{eqnarray*}

Taking $\varepsilon \leq \min \left \{ \dfrac{c_3}{2c_4}, \dfrac{c_0}{2n \delta_1} \right \}$ completes the proof. \qed

\hfill

Recall that a compact complex manifold $X$ is said to be \textit{Moishezon} if $X$ is bimeromorphic to a projective manifold.  Any compact complex manifold with a big line bundle $\mathcal{L} \to X$ is Moishezon, with the linear system $| \mathcal{L} |$ furnishing the bimeromorphic map.  In particular, a compact complex manifold $X$ of general type (i.e.,  the canonical bundle $K_X$ is big) is Moishezon.  In particular,  we have:

\subsection*{Corollary 2.7}
Let $(X, \omega)$ be a compact K\"ahler manifold with $(\delta_1, \delta_2)$--bounded geometry and $(\varepsilon, \delta)$--quasi-negative holomorphic sectional curvature. If one further assumes that  $K_X$ is nef,  then $X$ is of general type, and thus, in particular, Moishezon. \\

Moishezon's theorem \cite{Moishezon} asserts that a Moishezon manifold that is K\"ahler,  is projective. Further, from the recent developments in the minimal model program \cite[Corollary 1.4.6]{BCHM}, a Moishezon manifold without rational curves is projective.  In particular,  from \cite[Theorem 0.5]{DemaillyPaun}, we have the following:

\subsection*{Corollary 2.8}
Let $X$ be a compact K\"ahler manifold endowed with a Hermitian metric of $(\delta_1, \delta_2)$--bounded geometry and $(\varepsilon, \delta)$--quasi-negative holomorphic sectional curvature. If one further assumes that  $K_X$ is nef and $X$ is Kobayashi hyperbolic, then $X$ is projective with ample canonical bundle.

\subsection*{Remark 2.9}
It would be interesting to further explore the role played by Tian's $\alpha$--invariant, given its place in algebraic geometry.  Indeed,  we remind the reader that Tian's $\alpha$--invariant is an asymptotic version of the log canonical threshold (see \cite{TianAlpha, DK, CS}).

\section{Proof of Theorem 1.2}
To see \nameref{Theorem3}, we first prove a useful lemma by employing a method of Li--Ni--Zhu \cite{LiNiZhu}.  Let $$\lambda_{\omega_0} \ : = \ \max_{x \in X} \max_{v \in T_x^{1,0}X} \text{Ric}_{\alpha,\beta}^{\perp}(v).$$ Define \begin{eqnarray*}
\tau_{\omega_0}(x) & : = & \rho \left( \max_{v \in T_x^{1,0} X} \text{Ric}_{\alpha, \beta}^{\perp}(v) \right) \cdot \max_{v \in T_x^{1,0} X} \text{Ric}_{\alpha, \beta}^{\perp}(v),
\end{eqnarray*}

where $\rho : \mathbf{R} \to \{ n+1, 2n \}$ is a function with $\rho(s) = n+1$ for $s \leq 0$ and $\rho(s) = 2n$ for $s>0$.

\subsection*{Lemma 3.1}\label{key lem}
Let $(X^{n},\omega_{0})$ be a compact K\"{a}hler manifold. We consider the following complex Monge-Amp\`{e}re equation for $\varphi_t$:
\begin{equation}\label{MA equation}
(t\omega_{0}-\text{Ric}_{\omega_{0}}+ dd^c\varphi_t)^{n}
=e^{(1+\frac{\alpha}{2\beta})\varphi_t}\omega_{0}^{n}
\end{equation}
\begin{itemize}
\item[(a)] If $\lambda_{\omega_{0}}\geq0$, then the solution $\varphi_t$ to equation (\ref{MA equation}) exists for $t\in \left( \frac{2n\lambda_{\omega_{0}}}{(n+1)\alpha},+\infty \right)$;
\item[(b)] If $\lambda_{\omega_{0}}<0$, then the solution $\varphi_t$ to equation (\ref{MA equation}) exists for $t\in [0,+\infty)$.
\end{itemize}

\begin{proof} 
We consider the following complex Monge-Amp\`{e}re equation for $\varphi_t$:
\begin{equation}\label{MA equation1}
(t\omega_{0}-\text{Ric}_{\omega_{0}}+ dd^c \varphi_t)^{n}
=e^{(1+\frac{\alpha}{2\beta})\varphi_t}\omega_{0}^{n}.
\end{equation}
if $t$ is large enough, then $t\omega_{0}-\text{Ric}_{\omega_{0}}>0$. By the Aubin--Yau theorem \cite{Aubin,Yau1976}, there is  a unique smooth solution $\varphi_t$ of (\ref{MA equation1}) for $t\in (T,+\infty)$, where $T$ is the minimum existence time of $\varphi_t$.
We assume \nameref{key lem}(a) fails, i.e. $T>\frac{2n\lambda_{\omega_{0}}}{(n+1)\alpha}\geq 0$. Without losing generality, we only consider the $t\in (T,A]$, where $A$ is a fixed positive constant large enough.
Set $\omega_t : = t \omega_{0}- \text{Ric}_{\omega_{0}}+dd^c \varphi_t$.
Then for $t\in (T,A]$, we have
\begin{equation}\label{3.2}
\omega_t \ = \ t\omega_{0}-\text{Ric}_{\omega_{0}}+dd^c \varphi_t \ > \ 0,
\end{equation}
and (\ref{MA equation1}) is equivalent to
\begin{equation}\label{3.3}
\text{Ric}_{\omega_t} \ = \ -\omega_t +t\omega_{0}-
 \frac{\alpha}{2\beta} dd^c \varphi_t.
\end{equation}
To simplify notation we write
 the components of $\omega_t$ as $g_{i\overline{j}}$ and the components of $\omega_{0}$ as $h_{\alpha\overline{\beta}}$.
Let $G= \text{tr}_{\omega_t}\omega_{0}= \text{tr}_{g}h$,  from \nameref{SchwarzLemma} we have
\begin{equation}\label{3.4}
\Delta_{\omega_t}\log G \ \geq \ \frac{1}{G}(\text{Ric}_{i\overline{j}}^{g}
g^{i\overline{q}}g^{p\overline{j}}h_{p\overline{q}}-
g^{i\overline{j}}g^{k\overline{l}}R^{h}_{i\overline{j}k\overline{l}}).
\end{equation}
By \cite[Lemma 2.2]{LiNiZhu}, we get
\begin{equation}\label{3.5}
\frac{1}{G}g^{i\overline{j}}g^{k\overline{l}}R^{h}_{i\overline{j}k\overline{l}} \ \leq \ 
\frac{\lambda}{2\beta}G+\frac{\lambda}{2\beta G}|h|_{g}^{2}-\frac{\alpha}{\beta} \text{tr}_{g} \text{Ric}^{h}+\frac{\alpha}{2\beta G}(G\cdot \text{tr}_{g} \text{Ric}^{h}-\langle \omega_{h}, \text{Ric}^{h}\rangle_{g}).
\end{equation}
Choosing local coordinates such that $g_{i\overline{j}}=\delta_{i\overline{j}}$,
$h_{i\overline{j}}=h_{i\overline{i}}\delta_{ij}$, then
\begin{align}\label{3.6}
G\cdot \text{tr}_{g} \text{Ric}^{h}-\langle \omega_{h}, \text{Ric}^{h}\rangle_{g}& \ = \ \sum_{i} \text{Ric}_{i\overline{i}}^{h}
(\sum_{j}h_{j\overline{j}}-h_{i\overline{i}})=\sum_{i}(\text{Ric}^{h}_{i\overline{i}}
(\sum_{j\neq i}h_{j\overline{j}})) \nonumber \\
&\ \leq \  \sum_{i}(t\cdot h_{i\overline{i}}+\sqrt{-1}\varphi_{t,i\overline{i}})
(\sum_{j}h_{j\overline{j}}-h_{i\overline{i}})  \\
& \ = \ tG^{2}-t|h|_{g}^{2}+G\Delta_{g}\varphi_t-\langle dd^c \varphi_t, h\rangle_{g}. \nonumber
\end{align}
Here we used (\ref{3.2}). Combining (\ref{3.5}) and (\ref{3.6}), we get
\begin{align}\label{3.7}
\frac{1}{G}g^{i\overline{j}}g^{k\overline{l}}R^{h}_{i\overline{j}k\overline{l}} \ \leq \ 
&\frac{\lambda+\alpha t}{2\beta}G+\frac{\lambda-\alpha t}{2\beta}\frac{|h|_{g}^{2}}{G}-\frac{\alpha}{\beta} \text{tr}_{g} \text{Ric}^{h} \nonumber\\
&+\frac{\alpha}{2\beta}\Delta_{g}\varphi_t-\frac{\alpha}{2\beta}\frac{1}{G}\langle dd^c \varphi_t, h\rangle_{g}.
\end{align}
By using (\ref{3.3}), we have
\begin{align}\label{3.8}
\frac{1}{G} \text{Ric}_{i\overline{j}}^{g}
g^{i\overline{q}}g^{p\overline{j}}h_{p\overline{q}}&\ = \ \frac{1}{G}\langle
\text{Ric}^{g},h\rangle_{g} \ = \ \frac{1}{G}\langle-g+th-\frac{\alpha}{2\beta}
dd^c \varphi_t, h\rangle_{g} \nonumber \\
&\ = \ \frac{1}{G}\langle-g+th, h\rangle_{g}-\frac{1}{G}\frac{\alpha}{2\beta}\langle \varphi_t,  h\rangle_{g}.
\end{align}
Combining (\ref{3.4}),  (\ref{3.7}) and (\ref{3.8}),  we get
\begin{align}\label{3.9}
\Delta_{\omega_t} \left(\log G-\frac{\alpha}{2\beta}\varphi_t \right)  \ \geq \ &-\frac{\lambda+\alpha t}{2\beta}G-\frac{\lambda-\alpha t}{2\beta}\frac{|h|_{g}^{2}}{G}+\frac{\alpha}{\beta}(\text{tr}_{g}\text{Ric}^{h}-
\Delta_{g}\varphi_t) \nonumber\\
&+\frac{1}{G}\langle-g+th, h\rangle_{g}.
\end{align}
If $\lambda\geq 0$, since $t>T>0$, $n|h|_{g}^{2}\geq G^{2}$, $|h|_{g}^{2}\leq G^{2}$,
\begin{align}\label{3.10}
\Delta_{g}\varphi_t \ = \ n-tG+ \text{tr}_{g} \text{Ric}^{h},
\end{align}
and
\begin{align}\label{3.11}
\frac{1}{G}\langle-g+th, h\rangle_{g} \ = \ -1+\frac{t}{G}|h|_{g}^{2}\ \geq \ -1.
\end{align}
so we have
\begin{align}\label{3.12}
\Delta_{\omega_t} \left( \log G-\frac{\alpha}{2\beta}\varphi_t \right) \ \geq \ 
 \left( \frac{(n+1)\alpha t-2n\lambda}{2\beta n} \right) G-\left( \frac{\alpha n}{\beta}+1 \right).
\end{align}
Now we apply the maximum principle to get a uniform lower estimate of $\omega_t$.
For any fixed $t\in(T,A]$, at the maximum point $x_{1}$ of the $\varphi_t$,
we have $dd^c \varphi_t \leq 0$, so we can get
\begin{align}\label{3.13}
0 \ < \ g \ = \ \omega_t \mid_{x_{1}} \ \leq \ t\omega_{0}-\text{Ric}_{\omega_{0}}(x_{1}),
\end{align}
and
\begin{align}\label{3.14}
e^{(1+\frac{\alpha}{2\beta})\sup_{X}\varphi_t(x)} \ = \ e^{(1+\frac{\alpha}{2\beta})
\varphi_t(x_{1})} \ \leq \ \frac{(t\omega_{0}-\text{Ric}_{\omega_{0}}))^{n}}{\omega_{0}^{n}} \ \leq \  C,
\end{align}
where $C$ is independent of $t$. Hence $\varphi_t$ has a  uniform upper bound,  and we also get
\begin{align}\label{3.15}
\omega_t^{n} \ \leq \ C'\omega_{0}^{n}.
\end{align}
Again applying maximum principle to $\log G-\frac{\alpha}{2\beta}\varphi_t$.
At the point $x_{2}$, where the maximum of $\log G-\frac{\alpha}{2\beta}\varphi_t$ is attained, since $t>T>\frac{2n\lambda_{\omega_{0}}}{(n+1)\alpha}\geq 0 $, by (\ref{3.12}), we have that
\begin{align}\label{3.16}
G(x_{2}) \ \leq \ C'',
\end{align}
where $C''$ is also independent of $t$.
We note that
\begin{align}\label{3.17}
G\cdot \left( \frac{\omega_{0}^{n}}{\omega_t^{n}} \right)^{\frac{\alpha}{2\beta+\alpha}} \ = \ G\cdot \left(e^{-(1+\frac{\alpha}{2\beta})\varphi_t} \right)^{\frac{\alpha}{2\beta+\alpha}} \ = \ 
G\cdot e^{-\frac{\alpha}{2\beta}\varphi_t}.
\end{align}
So $\sup_{X} G\cdot \left( \frac{\omega_{0}^{n}}{\omega_t^{n}} \right) ^{\frac{\alpha}{2\beta+\alpha}}$
is also attained at $x_{2}$.  By the AM-GM inequality, at $x_{2}$, we have
$G\cdot \left( \frac{\omega_{0}^{n}}{\omega_t^{n}} \right)^{\frac{\alpha}{2\beta+\alpha}}\leq
G\cdot(\frac{G}{n})^{\frac{\alpha n}{2\beta+\alpha}}$. By (\ref{3.16}), we have
$\sup_{X}G\cdot(\frac{\omega_{0}^{n}}{\omega_t^{n}})^{\frac{\alpha}{2\beta+\alpha}}
\ \leq \ C'''$. Again by (\ref{3.15}), we get
\begin{align}\label{3.18}
\omega_t \ \geq \ \widetilde{C}\omega_{0},
\end{align}
where $\widetilde{C}$ is a positive number and independent of $t$. Combining (\ref{3.15}) and (\ref{3.18}), for some positive $C>0$ independent of $t$, we have
\begin{align}\label{3.19}
C^{-1} \omega_{0} \ \leq \ \omega_t \ \leq \ C\omega_{0}.
\end{align}
Therefore we obtain the higher order estimates (see \cite{TosattiYang}). Moreover, $\omega_t$ converges to a smooth K\"{a}hler form as $t\rightarrow T$. This
is a contradiction since $T$ is the minimum existence time of $\varphi_t$.  Therefore, there must holds $T\leq\frac{2n\lambda}{(n+1)\alpha}$.  This completes the proof of \nameref{key lem} (a).

If $\lambda<0$, by (\eqref{3.9}), it is easy to get that
\begin{align}\label{3.20}
\Delta_{\omega_t} \left( \log G-\frac{\alpha}{2\beta}\varphi_t \right) \ \geq \ 
\left( \frac{(n+1)\alpha t-(n+1)\lambda}{2\beta n} \right) G-\left( \frac{\alpha n}{\beta}+1\right).
\end{align}
Following the proof of  \nameref{key lem} (a), we can easily prove  \nameref{key lem} (b). 
\end{proof}

\noindent {\bf Proof of \nameref{Theorem3}:} The proof of  \nameref{Theorem3} is similar to  that of \nameref{Theorem2}. For completeness, we give a detailed proof.  We consider the following complex Monge-Amp\`{e}re equation for $\varphi_t$:
\begin{equation}\label{MA equation2}
(t\omega_{0}-\text{Ric}_{\omega_{0}}+ dd^c \varphi_t)^{n}
=e^{(1+\frac{\alpha}{2\beta})\varphi_t}\omega_{0}^{n}.
\end{equation}
let
$\omega_t =t\omega_{0}-\text{Ric}_{\omega_{0}}+ dd^c \varphi_t$.
By \nameref{key lem}, we have that
\begin{align}\label{3.22}
\Delta_{\omega_t} \left( \log (\text{tr}_{\omega_t} \omega_{0}) -\frac{\alpha}{2\beta}\varphi_t \right) \ \geq \ 
\left( \frac{(n+1)\alpha t-\tau_{\omega_{0}}}{2\beta n} \right) \text{tr}_{\omega_t}(\omega_{0})- \left( \frac{\alpha n}{\beta}+1 \right).
\end{align}
and there exists smooth solution $\varphi_t$ for $t\in \left( \frac{2n\lambda_{\omega_{0}}}{(n+1)\alpha},+\infty \right)$. In the following we just consider $t\in \left( \frac{2n\lambda_{\omega_{0}}}{(n+1)\alpha}, \frac{4n\lambda_{\omega_{0}}}{(n+1)\alpha} \right]$.  Integrating (\ref{3.22}) gives
\begin{align}\label{3.23}
\left( \frac{\alpha n}{\beta}+1 \right) \int_{X}\omega_t^{n}& \ \geq \  \int_{X}\left( \frac{(n+1)\alpha t}{2\beta n}-\frac{\tau_{\omega_{0}}}{2\beta n} \right)(\text{tr}_{\omega_t}\omega_{0})\omega_t^{n} \nonumber \\
& \ \geq \ -\frac{1}{2\beta n}\int_{\mathcal{U}}\tau_{\omega_{0}}(\text{tr}_{\omega_t}\omega_{0})\omega_t^{n} \nonumber \\
&\ \geq \ -\frac{1}{2\beta}\int_{\mathcal{U}}\tau_{\omega_{0}}\cdot \left( \frac{\omega_{0}^{n}}
{\omega_t^{n}} \right)^{\frac{1}{n}}\omega_t^{n} \nonumber \\
&=-\frac{1}{2\beta}\int_{\mathcal{U}}\tau_{\omega_{0}}\cdot e^{-\frac{1}{n}(1+\frac{\alpha}{2\beta})\varphi_t}\omega_t^{n} \nonumber \\
&\ \geq \ -\frac{1}{2\beta}e^{-\frac{1}{n}(1+\frac{\alpha}{2\beta})\sup_{X}\varphi_t}
\int_{\mathcal{U}}\tau_{\omega_{0}}\omega_t^{n}
\end{align}
Now let us estimate the bounds of the following fractions:
\begin{align}\label{3.24}
\frac{-\int_{\mathcal{U}}\tau_{\omega_{0}}\omega_t^{n}}{\int_{X}\omega_t^{n}} \ =\ \frac{-\int_{\mathcal{U}}\tau_{\omega_{0}}e^{(1+\frac{\alpha}{2\beta})\varphi_t}\omega_{0}^{n}}
{\int_{X}e^{(1+\frac{\alpha}{2\beta})\varphi_t}\omega_{0}^{n}}=
\frac{-\int_{\mathcal{U}}\tau_{\omega_{0}}e^{(1+\frac{\alpha}{2\beta})\varphi_t^{\ast}}\omega_{0}^{n}}
{\int_{X}e^{(1+\frac{\alpha}{2\beta})\varphi_t^{\ast}}\omega_{0}^{n}},
\end{align}
where $\varphi_t^{\ast} :=\varphi_t-\sup_{X}\varphi_t$. We fix a nonnegative number $b_{0} \geq 0$ such that $\text{Ric}_{\omega_{0}}\geq-b_{0}\omega_{0}$. Applying the maximum principle to $\varphi_t$, by (\ref{MA equation2}), we get that
\begin{align}\label{3.25}
\sup_{X}\varphi_t & \ \leq \ \frac{2\beta}{\alpha+2\beta}\log \left( \frac{4n\varepsilon}{(n+1)\alpha}+b_{0} \right)^{n} \nonumber \\
&\leq\frac{2\beta}{\alpha+2\beta}\log(c_{0}+b_{0})^{n}
\end{align}
for $t\in \left( \frac{2n\lambda_{\omega_{0}}}{(n+1)\alpha}, \frac{4n\lambda_{\omega_{0}}}{(n+1)\alpha} \right]$, where $c_{0}$ is a positive number
such that the $\alpha$--invariant $\alpha(X,c_{0}\omega_{0})\geq 2$ and we choose
$\varepsilon\leq\frac{(n+1)\alpha c_{0}}{4n}$.
Set
$$c_1 \ := \ \inf\left\{\int_{\mathcal{U}} e^{u}\omega_0^n|u\in \text{PSH}(X,(c_0+b_0)\omega_0),  \ \sup_Xu=0\right\},$$
which is a positive number depending only on $\mathcal{U},\omega_0,c_0,b_0,$, and hence only on $\mathcal{U},\omega_0$. Then, for $t\in \left( \frac{2n\lambda_{\omega_{0}}}{(n+1)\alpha}, \frac{4n\lambda_{\omega_{0}}}{(n+1)\alpha} \right]$, we have
\begin{align}
-\int_{\mathcal{U}}\tau_{\omega_{0}}e^{(1+\frac{\alpha}{2\beta})\varphi_t^{\ast}}\omega_{0}^{n} \ \geq \ 
(n+1)\delta\int_{\mathcal{U}}e^{(1+\frac{\alpha}{2\beta})
\varphi_t^{\ast}}\omega_{0}^{n}\ \geq \ (n+1)\delta c_{1} \nonumber
\end{align}
Since $\sup_{X}\varphi_t^{\ast}=0$, so we get 
\begin{align}
\int_{X}e^{(1+\frac{\alpha}{2\beta})\varphi_t^{\ast}}\omega_{0}^{n}\leq
\int_{X}\omega_{0}^{n} \ =: \ c_{2}, \nonumber
\end{align}
where $c_{2}$ is a positive number depending only on $\omega_{0}$. In conclusion,
we obtain that for $t\in \left( \frac{2n\lambda_{\omega_{0}}}{(n+1)\alpha}, \frac{4n\lambda_{\omega_{0}}}{(n+1)\alpha} \right]$,
\begin{align}
\frac{-\int_{\mathcal{U}}\tau_{\omega_{0}}\omega_t^{n}}{\int_{X}\omega_t^{n}}
\ \geq \ \frac{(n+1)\delta c_{1}}{c_{2}}. \nonumber
\end{align}
with, combining with (\ref{3.23}), concludes
\begin{align}
\sup_{\mathcal{U}}\varphi_t \ \geq  \ \frac{2\beta n}{\alpha+2\beta}\log\frac{(n+1)\delta c_{1}}{(2\alpha n+2\beta)c_{2}}.\nonumber
\end{align}
Now we define
\begin{align}
c_3\ : = \ \inf\left\{\int_Xe^v\omega_0^n \ | \ v\in \text{PSH}(X,c_0b_0\omega_0),  \ n\log\frac{(n+1)\delta c_{1}}{(2\alpha n+2\beta)c_{2}}\le\sup_Xv\le\log(c_0+b_0)^n\right\}, \nonumber
\end{align}
which is a positive number depending only on $U,\omega_0,\delta$.
Therefore,
\begin{align}\label{3.26}
\int_{X} \left( \frac{2n\lambda_{\omega_{0}}}{(n+1)\alpha}\omega_{0}-\text{Ric}_{\omega_{0}} \right)^{n}&=
\lim_{t\rightarrow \frac{2n\lambda_{\omega_{0}}}{(n+1)\alpha} }\int_{X} \left(t\omega_{0}-\text{Ric}_{\omega_{0}}+ dd^c \varphi_t \right)^{n} \nonumber \\
&=\lim_{t\rightarrow \frac{2n\lambda_{\omega_{0}}}{(n+1)\alpha}}\int_{X}
e^{(1+\frac{\alpha}{2\beta})\varphi_t}\omega_{0}^{n} \nonumber \\
&\geq \liminf_{t\rightarrow \frac{2n\lambda_{\omega_{0}}}{(n+1)\alpha}}\int_{X}
e^{(1+\frac{\alpha}{2\beta})\varphi_t}\omega_{0}^{n} \nonumber \\
&\geq c_{3}
\end{align}
On the other hand,
\begin{align}\label{3.27}
\int_{X} \left( \frac{2n\lambda_{\omega_{0}}}{(n+1)\alpha}\omega_{0}-\text{Ric}_{\omega_{0}} \right)^{n}&=
\lim_{t\rightarrow \frac{2n\lambda_{\omega_{0}}}{(n+1)\alpha} }\int_{X} \left( t\omega_{0}-\text{Ric}_{\omega_{0}}+dd^c \varphi_t \right)^{n} \nonumber \\
&\leq\int_{X} \left( \frac{4n\varepsilon}{(n+1)\alpha}\omega_{0}-\text{Ric}_{\omega_{0}} \right)^{n}\nonumber\\
&\leq \int_{X}(-\text{Ric}_{\omega_{0}})^{n}+c_{4}\varepsilon \nonumber\\
&=\int_{X}(2\pi c_{1}(K_{X}))^{n}+c_{4}\varepsilon
\end{align}
where $c_{4}$ is positive number depending only on $\omega_{0}$.
by (\ref{3.26}) and (\ref{3.27}), we have that
\begin{align}
\int_{X}(2\pi c_{1}(K_{X}))^{n}\geq c_{3}-c_{4}\varepsilon. \nonumber
\end{align}
Therefore, if $\varepsilon\leq \min\{\frac{c_{3}}{2c_{4}},\frac{(n+1)\alpha c_{0})}{4n}\}$, we have $\int_X(2\pi c_1(K_X))^n\ge\frac{c_3}{2}>0$.  \qed

\subsection*{Proof of \nameref{Theorem4}}
We note that the argument can be repeated more or less without change to prove \nameref{Theorem4}.  The only modification is that one requires \cite[Lemma 2.1]{LiNiZhu}.


\begin{thebibliography}{9}



\bibitem{Aubin}
Aubin, T., \'Equations du type Monge--Amp\`ere sur les variet\'es k\"ahl\'eriennes compactes, C.  Rendus Acad. Sci. Paris \textbf{283} (1976),  pp. 119–121.

\bibitem{BCHM}
Birkar,  C.,  Cascini, P.,  Hacon,  C.,  McKernan, J.,  Existence of minimal models for varieties of log general type. J. Amer. Math. Soc., 23(2):405--468, 2010.

\bibitem{BroderSBC}
Broder, K., {The Schwarz lemma in K\"ahler and non-K\"ahler geometry}, arXiv:2109.06331, to appear in the Asian J. Math.

\bibitem{BroderSBC2}
Broder, K.,  {The Schwarz lemma: An Odyssey}, arXiv:2110.04989,  Rocky Mountain J. Math. 52 (2022), no. 4, 1141–1155.

\bibitem{BroderRemarks}
Broder, K.,  Remarks on the Wu--Yau theorem,  (submitted).

\bibitem{BroderStanfield}
Broder, K., Stanfield, J.,  On the Gauduchon Curvature of Hermitian Manifolds,  Int. J.  of Math., doi: \url{https://doi.org/10.1142/S0129167X23500398}

\bibitem{BT}
Broder, K.,  Tang, K., On the weighted orthogonal Ricci curvature,  J. Geom. and Phys.,  104783, doi: \url{https://doi.org/10.1016/j.geomphys.2023.104783}.

\bibitem{BroderTangAltered}
Broder, K., Tang, K., On the altered holomorphic curvatures of Hermitian manifolds,  arXiv:2201.03666.

\bibitem{Campana}
Campana, F.,  Twistor spaces and nonhyperbolicity of certain symplectic K\"ahler manifolds, In Complex analysis (Wuppertal 1991), Aspects Math., E17, pp. 64--69. Friedr. Veiweg, Braunschweig, 1991.

\bibitem{CS}
Cheltsov, I., Shramov, C., Log canonical thresholds of smooth Fano threefolds (Russian), with an appendix by Jean-Pierre Demailly, Uspekhi Mat. Nauk \textbf{63} (2008), pp. 73--180; translation in Russian Math. Surveys \textbf{63} (2008), pp. 859--958.

\bibitem{ChuLeeTam}
Chu, J. C.,  Lee,  M. C., Tam L.-F., K\"ahler manifolds and mixed curvature. Trans. Amer. Math. Soc. 375 (2022), no. 11, 7925-7944.

\bibitem{Debarre}
Debarre,  O., Higher-dimensional algebraic geometry, Universitext, Springer-Verlag, New York, 2001. MR 1841091



\bibitem{DK}
Demailly, J.-P., Koll\'ar, J.,  Semi-continuity of complex singularity exponents and K\"ahler--Einstein metrics on Fano orbifolds, Ann. Sci. \'Ecole Norm.  Sup. (4) \textbf{34} (2001),  525--556.

\bibitem{DemaillyPaun}
Demailly, J.-P.,  Paun, M., Numerical characterization of the K\"ahler cone of a compact K\"ahler manifold, Ann. of Math. (2) \textbf{159} (2004), no. 3, 1247--1274.



\bibitem{DiverioTrapani}
Diverio, S.,  Trapani, S., Quasi-negative holomorphic sectional curvature and positivity of the canonical bundle, J. Diff. Geom., \textbf{111} (2019), pp. 303--314.



\bibitem{HeierLuWong}
Heier, G., Lu, S., Wong, B., On the canonical line bundle and negative holomorphic sectional curvature,  Math. Res. Lett.  \textbf{17} (2010),  no. 06,  pp. 1101–1110

\bibitem{LeeStreets}
Lee, M. C., Streets, J., Complex manifolds with negative curvature operator, Int. Math. Res. Not,  arXiv: 1903.12645.

\bibitem{LiNiZhu}
Li,  C., Ni, L., Zhu, X., An application of a $C^2$-estimate for a complex Monge-Amp\`ere equation,  Internat.  J.  Math.  32 (2021), no. 12, Paper No. 2140007, 13 pp. 

\bibitem{Moishezon}
Moishezon,  B. G., On $n$--dimensional compact complex manifolds having $n$ algebraically independent meromorphic functions. I. Izv. Akad. Nauk SSSR Ser. Mat., 30: 133--174, 1966.

\bibitem{NikRicci}
Ni, L.,  Liouville theorems and a Schwarz Lemma for holomorphic mappings between K\"ahler manifolds.  Comm.  Pure Appl.  Math, 2021, 74: 1100--1126.

\bibitem{Royden}
Royden, H. L.,  The Ahlfors-Schwarz lemma in several complex variables, Comment. Math. Helvetici \textbf{55} (1980), p. 547--558.

\bibitem{SiuMoishezon}
Siu, Y.-T.,  A vanishing theorem for semipositive line bundles over non-K\"ahler manifolds, J. Diff. Geom. \textbf{19} (1984), pp. 431--452.

\bibitem{TangkRicci}
Tang, K., On almost nonpositive $k$-Ricci curvature,  J. Geom. Anal. 32 (2022), no. 12, Paper No. 306, 12 pp.

\bibitem{TianAlpha}
Tian, G.,  On K\"ahler-Einstein metrics on certain K\"ahler manifolds with $C_1(M) > 0$, Invent.  Math. 89 (1987), no. 2, 225--246

\bibitem{TosattiYang}
Tosatti, V., Yang, X.-K., An extension of a theorem of Wu--Yau, J. Differential Geom. 107(3): 573--579

\bibitem{Wong}
Wong, B., The uniformization of compact K\"ahler surfaces of negative curvature, J. Differential Geom. 16 (1981), no. 3, 407--420. 

\bibitem{WuYau1}
Wu, D., Yau, S.-T., Negative holomorphic sectional curvature and positive canonical bundle, Invent. Math. \textbf{204} (2016), no. 2, 595--604

\bibitem{WuYau2}
Wu, D., Yau, S.-T., A remark on our paper "Negative holomorphic curvature and positive canonical bundle". Comm Anal Geom. \textbf{24} (2016) 901--912.

\bibitem{YangZhengCurvature}
Yang, B., Zheng, F., {On curvature tensors of Hermitian manifolds}, Communications in Analysis and Geometry, vol. \textbf{26} (2018), no. 5, pp. 1195--1222

\bibitem{YangZhengRBC}
Yang, X., Zheng, F., On the real bisectional curvature for Hermitian manifolds, Trans. Amer. Math. Soc. \textbf{371} (2019), no. 4, 2703--2718

\bibitem{Yau1976}
Yau, S.-T., On the Ricci curvature of a compact K\"ahler manifold and the complex Monge--Amp\`ere equation, I, Comm. Pure Appl. Math. 31 (1978),  pp. 339–411.

\bibitem{Zhang}
Zhang, Y.,  Integral inequalities for holomorphic maps and applications, Trans. Amer. Math. Soc. 374 (2021), no.4, 2341--2358

\bibitem{Zhang1}
Zhang, Y., Holomorphic sectional curvature, nefness and Miyaoka-Yau type inequality, Math. Z. 298 (2021), no.3-4, 953–974

\bibitem{ZhangZheng}
Zhang, Y., Zheng, T.,  On almost quasi-negative holomorphic sectional curvature,  arXiv:2010.01314v4.
\end{thebibliography}
\end{document}